\documentclass[amstex,12pt, amssymb]{article}

\usepackage{mathtext}
\usepackage[cp1251]{inputenc}
\usepackage[T2A]{fontenc}
\usepackage[dvips]{graphicx}
\usepackage{amsmath}
\usepackage{amssymb}
\usepackage{amsxtra}
\usepackage{latexsym}
\usepackage{ifthen}

\newcounter{lemma}[section]

\newcounter{corollary}[section]

\newcounter{remark}[section]

\newcounter{theorem}[section]

\newcounter{proposition}[section]

\newcounter{example}

\numberwithin{equation}{section}

\begin{document}

\markboth{E.~SEVOST'YANOV, V. DESYATKA, Z.~KOVBA}{\centerline{ON THE
PRIME ENDS EXTENSION ...}}

\def\cc{\setcounter{equation}{0}
\setcounter{figure}{0}\setcounter{table}{0}}

\overfullrule=0pt


\author{EVGENY SEVOST'YANOV, VICTORIA DESYATKA, ZARINA KOVBA}

\title{
{\bf ON THE PRIME ENDS EXTENSION OF UNCLOSED INVERSE MAPPINGS}}

\date{\today}
\maketitle

\begin{abstract}
We consider mappings that distort the modulus of families of paths
in the opposite direction in the manner of Poletsky's inequality.
Here we study the case when the mappings are not closed, in
particular, they do not preserve the boundary of the domain under
the mapping. Under certain conditions, we obtain results on the
continuous boundary extension of such mappings in the sense of prime
ends. In addition, we obtain corresponding results on the
equicontinuity of families of such mappings in terms of prime ends.
\end{abstract}

\bigskip
{\bf 2010 Mathematics Subject Classification: Primary 30C65;
Secondary 31A15, 31B25}

\section{Introduction}

This paper is devoted to the boundary behavior of mappings in terms
of prime ends. Note that this problem has been the subject of a
fairly large number of publications, including papers on
quasiconformal mappings and mappings with finite distortion, see,
for example, \cite{Ad, ABBS, Car, GU, GRY, KR} and~\cite{KrPa}. In
particular, the following statement holds
(see~\cite[Theorem~4.1]{Na}).

\medskip
{\bf Theorem.} {\sl Under a quasiconformal mapping $f$ of a collared
domain $D_0$ onto a domain $D,$ there exists a one-to-one
correspondence between the boundary points of $D_0$ and the prime
ends of $D.$ Moreover, the cluster set $C(f, b),$ $b\in \partial
D_0,$ coincides with the impression $I(P)$ of the corresponding
prime end $P$ of $D.$}

\medskip
Relatively recently, the first co-author obtained results on
continuous extension of mappings with branching that satisfy
Poletsky's inequality or its inverse, see, for example,
\cite{Sev$_2$}, \cite{Sev$_3$} and \cite{Sev$_1$}. Note that the
closeness of mappings is a very strong additional restriction,
significantly narrowing the class under consideration. The goal of
this article is to consider mappings which are nor closed. We show
that open discrete mappings have a continuous boundary extension in
terms of prime ends, although mappings may, in general, not be
closed. However, there are a number of additional restrictions both
on the domain under consideration and on the mappings themselves. We
also note two recent preprints by the first two co-authors, dealing
with similar mappings for the case of ``simple'' boundaries (here
the usual Euclidean continuous extension of mappings holds), see,
e.g., \cite{DS$_1$}--\cite{DS$_2$}.

\medskip
Let us recall the definition. Let $y_0\in {\Bbb R}^n,$
$0<r_1<r_2<\infty$ and
\begin{equation}\label{eq1**}
A=A(y_0, r_1,r_2)=\left\{ y\,\in\,{\Bbb R}^n:
r_1<|y-y_0|<r_2\right\}\,.\end{equation}
Given sets $E,$ $F\subset\overline{{\Bbb R}^n}$ and a domain
$D\subset {\Bbb R}^n$ we denote by $\Gamma(E,F,D)$ a family of all
paths $\gamma:[a,b]\rightarrow \overline{{\Bbb R}^n}$ such that
$\gamma(a)\in E,\gamma(b)\in\,F $ and $\gamma(t)\in D$ for $t \in
(a, b).$ If $f:D\rightarrow {\Bbb R}^n,$ $y_0\in f( D)$ and
$0<r_1<r_2<d_0=\sup\limits_{y\in f(D)}|y-y_0|,$ then by
$\Gamma_f(y_0, r_1, r_2)$ we denote the family of all paths $\gamma$
in $D$ such that $f(\gamma)\in \Gamma(S(y_0, r_1), S( y_0, r_2),
A(y_0,r_1,r_2)).$ Let $Q:{\Bbb R}^n\rightarrow [0, \infty]$ be a
Lebesgue measurable function. We say that {\it $f$ satisfies
Poletsky inverse inequality} at the point $y_0\in f(D),$ if the
relation
\begin{equation}\label{eq2*A}
M(\Gamma_f(y_0, r_1, r_2))\leqslant \int\limits_{A(y_0,r_1,r_2)\cap
f(D)} Q(y)\cdot \eta^n (|y-y_0|)\, dm(y)
\end{equation}
holds for any Lebesgue measurable function $\eta:
(r_1,r_2)\rightarrow [0,\infty ]$ such that
\begin{equation}\label{eqA2}
\int\limits_{r_1}^{r_2}\eta(r)\, dr\geqslant 1\,.
\end{equation}
Many mappings, including quasiregular mappings and mappings with
finite distortion, satisfy estimates~(\ref{eq2*A})--(\ref{eqA2}),
see, for example,~\cite[Theorem~1]{Pol}, \cite[Remark~2.5.II]{Ri}
and \cite[Theorems~8.1, 8.5]{MRSY}, cf.~\cite[Theorem~3]{Cr}.

\medskip
Recall that a mapping $f:D\rightarrow {\Bbb R}^n$ is called {\it
discrete} if the pre-image $\{f^{-1}\left(y\right)\}$ of each point
$y\,\in\,{\Bbb R}^n$ consists of isolated points, and {\it is open}
if the image of any open set $U\subset D$ is an open set in ${\Bbb
R}^n.$ Later, in the extended space $\overline{{{\Bbb R}}^n}={{\Bbb
R}}^n\cup\{\infty\}$ we use the {\it spherical (chordal) metric}
$h(x,y)=|\pi(x)-\pi(y)|,$ where $\pi$ is a stereographic projection
$\overline{{{\Bbb R}}^n}$ onto the sphere
$S^n(\frac{1}{2}e_{n+1},\frac{1}{2})$ in ${{\Bbb R}}^{n+1},$ namely,
$$h(x,\infty)=\frac{1}{\sqrt{1+{|x|}^2}}\,,$$
\begin{equation}\label{eq3C}
h(x,y)=\frac{|x-y|}{\sqrt{1+{|x|}^2} \sqrt{1+{|y|}^2}}\,, \quad x\ne
\infty\ne y
\end{equation}
(see \cite[Definition~12.1]{Va}). Further, the closure
$\overline{A}$ and the boundary $\partial A$ of the set $A\subset
\overline{{\Bbb R}^n}$ we understand relative to the chordal metric
$h$ in $\overline{{\Bbb R}^n}.$

\medskip
The boundary of $D$ is called {\it weakly flat} at the point $x_0\in
\partial D,$ if for every $P>0$ and for any neighborhood $U$
of the point $x_0$ there is a neighborhood $V\subset U$ of the same
point such that $M(\Gamma(E, F, D))>P$ for any continua $E, F\subset
D$ such that $E\cap\partial U\ne\varnothing\ne E\cap\ partial V$ and
$F\cap\partial U\ne\varnothing\ne F\cap\partial V.$ The boundary of
$D$ is called weakly flat if the corresponding property holds at any
point of the boundary $D.$

\medskip
Given a mapping $f:D\rightarrow {\Bbb R}^n$, we denote
\begin{equation}\label{eq1_A_4} C(f, x):=\{y\in \overline{{\Bbb
R}^n}:\exists\,x_k\in D: x_k\rightarrow x, f(x_k) \rightarrow y,
k\rightarrow\infty\}
\end{equation}
and
\begin{equation}\label{eq1_A_5} C(f, \partial
D)=\bigcup\limits_{x\in \partial D}C(f, x)\,.
\end{equation}
In what follows, ${\rm Int\,}A$ denotes the set of inner points of
the set $A\subset \overline{{\Bbb R}^n}.$ Recall that the set
$U\subset\overline{{\Bbb R}^n}$ is neighborhood of the point $z_0,$
if $z_0\in {\rm Int\,}A.$

\medskip
Consider the following definition, which goes back to
N\"akki~\cite{Na}, see also~\cite{KR}. The boundary of a domain $D$
in ${\Bbb R}^n$ is said to be {\it locally quasiconformal} if every
$x_0\in\partial D$ has a neighborhood $U$ that admits a
quasiconformal mapping $\varphi$ onto the unit ball ${\Bbb
B}^n\subset{\Bbb R}^n$ such that $\varphi(\partial D\cap U)$ is the
intersection of ${\Bbb B}^n$ and a coordinate hyperplane. The
sequence of cuts $\sigma_m,$ $m=1,2,\ldots ,$ is called {\it
regular,} if
$\overline{\sigma_m}\cap\overline{\sigma_{m+1}}=\varnothing$ for
$m\in {\Bbb N}$ and, in addition, $d(\sigma_{m})\rightarrow 0$ as
$m\rightarrow\infty.$ If the end $K$ contains at least one regular
chain, then $K$ will be called {\it regular}. We say that a bounded
domain $D$ in ${\Bbb R}^n$ is {\it regular}, if $D$ can be
quasiconformally mapped to a domain with a locally quasiconformal
boundary whose closure is a compact in ${\Bbb R}^n,$ and, besides
that, every prime end in $D$ is regular. Note that space
$\overline{D}_P=D\cup E_D$ is metric, which can be demonstrated as
follows. If $g:D_0\rightarrow D$ is a quasiconformal mapping of a
domain $D_0$ with a locally quasiconformal boundary onto some domain
$D,$ then for $x, y\in \overline{D}_P$ we put:
\begin{equation}\label{eq5M}
\rho(x, y):=|g^{\,-1}(x)-g^{\,-1}(y)|\,,
\end{equation}
where the element $g^{\,-1}(x),$ $x\in E_D,$ is to be understood as
some (single) boundary point of the domain $D_0.$ The specified
boundary point is unique and well-defined by~\cite[Theorem~2.1,
Remark~2.1]{IS}, cf.~\cite[Theorem~4.1]{Na}. It is easy to verify
that~$\rho$ in~(\ref{eq5M}) is a metric on $\overline{D}_P.$

\medskip
Assume that, $D^{\,\prime}\setminus C(f, \partial
D)=\bigcup\limits_{i=1}^ND_i,$ where $D_i$ is a regular domain for
$1\leqslant i\leqslant N,$ $D_i\cap D_j=\varnothing$ for $i\ne j.$
Let $\rho_i$ is a metric in $\overline{D_i}_P$ defined
by~(\ref{eq5M}). We set $\rho_j(x, y)=0$ for $x, y\in
\bigcup\limits_{i=1}^N\overline{D_i}_P\setminus \overline{D_j}_P,$
and $\rho_j(x, y)=1$ for $x\in \overline{D_j}_P$ and $y\not\in
\overline{D_j}_P.$ Then $\rho_j$ is a pseudometric on
$\bigcup\limits_{i=1}^N\overline{D_i}_P$ (see
\cite[2.2.21.XV]{Ku$_1$}). Set
\begin{equation}\label{eq1}
\rho(x, y)=\sum\limits_{i=1}^n2^{\,-i}\frac{\rho_i(x,
y)}{1+\rho_i(x, y)}\,.
\end{equation}
Observe that, the metric $\rho^{\,*}_j(x, y):=\frac{\rho_j(x,
y)}{1+\rho_j(x, y)}$ is a metric on $\overline{D_j}_P$ for
$1\leqslant j\leqslant N$ (see \cite[21.2.V]{Ku$_1$}) which is
homeomorphic to $\rho_j(x, y)$ on $\overline{D_i}_P$ (see
\cite[21.2.V]{Ku$_1$}). Besides that, $\rho(x, y)$ is a metric on
$\bigcup\limits_{i=1}^N\overline{D_i}_P,$ see
\cite[Remark~2.2.21.XV]{Ku$_1$}.

\medskip
Observe that, the space $\bigcup\limits_{i=1}^N\overline{D_i}_P$ is
compact. Indeed, if $x_k\in \bigcup\limits_{i=1}^N\overline{D_i}_P,$
then there is $1\leqslant j_0\leqslant N$ and a subsequence
$x_{k_l},$ $l=1,2,\ldots ,$ such that $x_{k_l}\in
\overline{D_{j_0}}_P$ for any $l\in {\Bbb N}.$ We may find a
sequence $y_l\in D_{j_0},$ $l=1,2,\ldots ,$ such that
$\rho_{j_0}(y_l, x_{k_l})<1/l$ for $l\in {\Bbb N}.$ Since by the
definition the domain $D_{j_0}$ is bounded, we may consider that
$y_l\rightarrow y_0\in \overline{D_{j_0}}.$ There are two cases:
$y_0\in D_{j_0},$ or $y_0\in \partial D_{j_0}.$ In the first case,
when $y_0\in D_{j_0},$ by the definition of $\rho_{j_0}$
in~(\ref{eq5M}), we have that $\rho_{j_0}(y_l, y_0)\rightarrow 0$ as
$l\rightarrow \infty.$ Now, by the triangle inequality
$\rho_{j_0}(x_{k_l}, y_0)\leqslant \rho_{j_0}(x_{k_l},
y_l)+\rho_{j_0}(y_l, y_0)<\frac{1}{l}+\rho_{j_0}(y_l,
y_0)\rightarrow 0$ as $l\rightarrow\infty.$ Thus,
$x_{k_l}\rightarrow y_0$ as $l\rightarrow \infty$ in the
metric~$\rho_{j_0}.$ By the definition of the metric $\rho$
in~(\ref{eq1}) $x_{k_l}\rightarrow y_0$ as $l\rightarrow \infty$ in
the metric~$\rho,$ as well. In the second case, when $y_0\in
\partial D_{j_0},$ we consider the homeomorphism $g:D_0\rightarrow
D_{j_0}$ of a domain $D_0$ with a locally quasiconformal boundary
onto $D_{j_0}.$ Such a domain $D_0$ and a homeomorphism $g$ exist by
the definition of a regular domain. Since $D_0$ is bounded, we may
consider that $g^{\,-1}(x_{k_l})\rightarrow z_0\in \partial D_0$ as
$l\rightarrow\infty.$ Set $P_0=g(z_0)\in E_{D_{j_0}}.$ The prime end
$P_0$ is well-defined due to~\cite[Theorem~2.1, Remark~2.1]{IS},
cf.~\cite[Theorem~4.1]{Na}. In addition, by the definition of the
metrics $\rho_{j_0},$ $\rho_{j_0}(x_{k_l}, P_0)\rightarrow 0$ as
$l\rightarrow\infty.$ By the definition of $\rho$ in~(\ref{eq1}),
$x_{k_l}\rightarrow P_0$ as $l\rightarrow \infty$ by the
metric~$\rho,$ as well.

\medskip
The following statement holds.

\medskip
\begin{theorem}\label{th3}
{\it\, Let $D$ and $D^{\,\prime}$ be domains in ${\Bbb R}^n,$
$n\geqslant 2,$ and let $D$ be a domain with a weakly flat boundary.
Suppose that $f$ is open discrete mapping of $D$ onto $D^{\,\prime}$
satisfying the relation~(\ref{eq2*A}) at each point $y_0\in
D^{\,\prime}.$ In addition, assume that the following conditions are
fulfilled:

\medskip
1) for each point $y_0\in \partial D^{\,\prime}$ there is
$0<r_0:=\sup\limits_{y\in D^{\,\prime}}|y-y_0|$ such that for any
$0<r_1<r_2<r_0:=\sup\limits_{y\in D^{\,\prime}}|y-y_0|$ there exists
a set $E\subset[r_1, r_2]$ of positive linear Lebesgue measure such
that $Q$ is integrable on $S(y_0, r)$ for $r\in E;$

\medskip
2) $D^{\,\prime}\setminus C(f,
\partial D)$ consists of finite number of pairwise disjoint
domains $D_1, D_2,\ldots, D_N$ any of which is regular;

\medskip
3) the set $f^{\,-1}(C(f, \partial D)\cap D^{\,\prime})$ is nowhere
dense in $D.$

\medskip
Then the mapping $f$ has a boundary extension
$\overline{f}:\overline{D}\rightarrow
\bigcup\limits_{i=1}^N\overline{D_i}_P$ which is continuous in
$D\setminus f^{\,-1}(C(f, \partial D)\cap D^{\,\prime})$ by the
metric $\rho$ defined in~(\ref{eq1}). Moreover,
$\overline{f}(\overline{D})=\bigcup\limits_{i=1}^N\overline{D_i}_P.$
}
\end{theorem}

\medskip
\begin{corollary}\label{cor1}
{\it\, The statement of Theorem~\ref{th3} remains true if, instead
of condition~1), a simpler condition holds: $Q\in
L^1(D^{\,\prime}).$}
\end{corollary}

\section{Proof of Theorem~\ref{th3}}

Let $D\subset {\Bbb R}^n,$ $f:D\rightarrow {\Bbb R}^n$ be a discrete
open mapping, $\beta: [a,\,b)\rightarrow {\Bbb R}^n$ be a path, and
$x\in\,f^{\,-1}(\beta(a)).$ A path $\alpha: [a,\,c)\rightarrow D$ is
called a {\it maximal $f$-lifting} of $\beta$ starting at $x,$ if
$(1)\quad \alpha(a)=x\,;$ $(2)\quad f\circ\alpha=\beta|_{[a,\,c)};$
$(3)$\quad for $c<c^{\prime}\leqslant b,$ there is no a path
$\alpha^{\prime}: [a,\,c^{\prime})\rightarrow D$ such that
$\alpha=\alpha^{\prime}|_{[a,\,c)}$ and $f\circ
\alpha^{\,\prime}=\beta|_{[a,\,c^{\prime})}.$ Here and in the
following we say that a path $\beta:[a, b)\rightarrow
\overline{{\Bbb R}^n}$ converges to the set $C\subset
\overline{{\Bbb R}^n}$ as $t\rightarrow b,$ if $h(\beta(t),
C)=\sup\limits_{x\in C}h(\beta(t), C)\rightarrow 0$ at $t\rightarrow
b.$ The following is true (see~\cite[Lemma~3.12]{MRV}).

\medskip
\begin{proposition}\label{pr3}
{\it\, Let $f:D\rightarrow {\Bbb R}^n,$ $n\geqslant 2,$ be an open
discrete mapping, let $x_0\in D,$ and let $\beta: [a,\,b)\rightarrow
{\Bbb R}^n$ be a path such that $\beta(a)=f(x_0)$ and such that
either $\lim\limits_{t\rightarrow b}\beta(t)$ exists, or
$\beta(t)\rightarrow \partial f(D)$ as $t\rightarrow b.$ Then
$\beta$ has a maximal $f$-lifting $\alpha: [a,\,c)\rightarrow D$
starting at $x_0.$ If $\alpha(t)\rightarrow x_1\in D$ as
$t\rightarrow c,$ then $c=b$ and $f(x_1)=\lim\limits_{t\rightarrow
b}\beta(t).$ Otherwise $\alpha(t)\rightarrow \partial D$ as
$t\rightarrow c.$}
\end{proposition}

\medskip
{\it Proof of Theorem~\ref{th3}.} We fix $x_0\in\partial D.$ Let us
show that the mapping $f$ has a continuous extension to a point
$x_0.$ Using the M\"{o}bius transformation $\varphi:\infty\mapsto 0$
if necessary, due to the invariance of the modulus $M$
in~(\ref{eq2*A}) (see~\cite[Theorem~8.1]{Va}), we may assume that
$x_0\ne\infty.$

\medskip
Fix $x_0\in \partial D.$ Let $x_i\in D,$ $i=1,2,\ldots,$ be a
sequence converging to $x_0$ as $i\rightarrow\infty.$ Since the set
$f^{\,-1}(C(f, \partial D)\cap D^{\,\prime})$ is nowhere dense in
$D,$ there exists a sequence $x_{ki}\in D\setminus(f^{\,-1}(C(f,
\partial D)\cap D^{\,\prime})),$
$k=1,2,\ldots ,$ such that $x_{ki}\rightarrow x_i$ as
$k\rightarrow\infty.$ Thus, for any $i\in {\Bbb N}$ there exists
$k_i\in {\Bbb N}$ such that $|x_{k_i}i-x_i|<2^{\,-i}.$ Consequently,
by the triangle inequality
$$|x_{{k_i}i}-x_0|\leqslant |x_{{k_i}i}-x_i|+|x_i-x_0|<2^{\,-i}+|x_i-x_0|
\rightarrow 0$$
as $i\rightarrow\infty.$ By the definition, $f(x_{{k_i}i})\in
\bigcup\limits_{i=1}^ND_i$ for each $i\in {\Bbb N}.$ Since the space
$\bigcup\limits_{i=1}^N\overline{D_i}_P$ is compact (see remarks
made before Theorem~\ref{th1}), we may consider that the sequence
$f(x_{{k_i}i})$ converges to some $P_0\in
\bigcup\limits_{i=1}^N\overline{D_i}_P$ as $i\rightarrow\infty.$ Set
\begin{equation}\label{eq1A}
f(x_0):=P_0\,.
\end{equation}
Now, let us to prove that the continuous extension in~(\ref{eq1A})
is continuous in $D\setminus f^{\,-1}(C(f, \partial D)\cap
D^{\,\prime})$ by the metric $\rho$ defined in~(\ref{eq1}).

\medskip
Let us carry out the proof of this fact from the opposite, namely,
suppose that $f$ has no a continuous extension to the point $x_0$
with respect to $D\setminus f^{\,-1}(C(f, \partial D)\cap
D^{\,\prime}).$ Due to the compactness of the space
$\bigcup\limits_{i=1}^N\overline{D_i}_P,$ there exist sequences
$x_i, y_i\in D\setminus f^{\,-1}(C(f, \partial D)\cap
D^{\,\prime}),$ $i=1,2,\ldots ,$ such that $x_i, y_i\rightarrow x_0$
as $i\rightarrow\infty,$ and
\begin{equation}\label{eq1B}
\rho(f(x_i), f(y_i))\geqslant a>0
\end{equation}
for some $a>0$ and all $i\in {\Bbb N},$ where $\rho$ is a metric
in~(\ref{eq1}). Again, by the compactness of
$\bigcup\limits_{i=1}^N\overline{D_i}_P,$ we may assume that the
sequences $f(x_i)$ and $f(y_i)$ converge as $i\rightarrow\infty$ to
$P_1$ and $P_2\in \bigcup\limits_{i=1}^N\overline{D_i}_P,$
respectively.
\medskip
By the assumption~2), $D^{\,\prime}\setminus C(f,
\partial D))$ consists of finite number of domains $D_1, D_2,\ldots, D_N$
any of which is regular. Consequently, there are domains $D_{j_1}$
and $D_{j_2,}$ $1\leqslant j_1, j_2\leqslant N,$ and subsequences
$f(x_{i_k}), f(y_{j_k}),$ $k=1,2,\ldots,$ of $f(x_i)$ and $f(y_i),$
such that $f(x_{i_k})\in D_{j_1}$ and $f(y_{i_k})\in D_{j_2}$ for
$k\in {\Bbb N}.$ Without loss of generality, making a relabeling if
necessary, we may assume that $f(x_i)\in D_{j_1}$ and $f(y_i)\in
D_{j_2},$ $i=1,2,\ldots .$ Here the case $j_1=j_2$ is possible.

\medskip
Since the spaces $\overline{D_{j_1}}_P$ and $\overline{D_{j_2}}_P$
is compact (see remarks made before the formulation of
Theorem~\ref{th3}), we may consider that $f(x_i)$ and $f(y_i)$
converge to $P_1\in \overline{D_{j_1}}_P$ and $P_2 \in
\overline{D_{j_2}}_P$ by the metrics $\rho_{j_1}$ and $\rho_{j_2},$
correspondingly. Observe that, by the definition of $C(f, D),$ that
$P_1\in E_{D_{j_1}}$ and $P_2\in E_{D_{j_2}}.$

\medskip
There are two cases: $j_1\ne j_2$ and $j_1=j_2.$ If $j_1\ne j_2,$
then obviously $P_1\ne P_2.$ Otherwise, if $j_1=j_2,$ then by the
definition of $\rho$ in (\ref{eq1}), $\rho(x_i,
y_i)=2^{\,-j_1}\frac{\rho_{j_1}(x_i, y_i)}{1+\rho_{j_1}(x_i, y_i)}.$
However, due to~(\ref{eq1B}), by the triangle inequality,
$$\rho(P_1, P_2)\geqslant \rho(P_2, f(x_i))-\rho(f(x_i), P_1)\geqslant
\rho(f(x_i), f(y_i))-\rho(f(y_i), P_2)-\rho(f(x_i), P_1)=$$
$$=\rho(f(x_i), f(y_i))- 2^{\,-j_1}\frac{\rho_{j_1}(f(y_i),
P_2)}{1+\rho_{j_1}(f(y_i), P_2)}-2^{\,-j_1}\frac{\rho_{j_1}(f(x_i),
P_1)}{1+\rho_{j_1}(f(x_i), P_1)}\geqslant a/2>0$$
for sufficiently large $i\in {\Bbb N}.$ Thus, $P_1\ne P_2,$ as well.
In any of two cases, $j_1\ne j_2$ or $j_1=j_2,$ we have that $P_1\ne
P_2.$

\medskip
Let $\sigma_i$ and $\sigma^{\,\prime}_i,$ $i=0,1,2,\ldots, $ be
sequences of cuts corresponding $P_1$ and $P_2.$ We may consider
that the cuts $\sigma_i,$ $i=0,1,2,\ldots, $ lie on the spheres
$S(z_0, r_i)$ centered at some point $z_0\in
\partial D_{j_1},$ where $r_i\rightarrow 0$ as $i\rightarrow\infty$
(such a sequence $\sigma_i$ exists due to~\cite[Lemma~3.1]{IS},
cf.~\cite[Lemma~1]{KR}). Let $d_i$ and $g_i,$ $i=0,1,2,\ldots, $ are
sequences of nested domains in $D_{j_1},$ corresponding to cuts
$\sigma_i$ and $\sigma^{\,\prime}_i.$ Since the space
$(\overline{D_{j_1}}_P, \rho_{j_1})$ is metric, we may consider that
$d_i$ and $g_i$ are disjoint, so that
\begin{equation}\label{eq4}
d_0\cap g_0=\varnothing\,.
\end{equation}
Since $f(x_i)$ converges to $P_1$ as $i\rightarrow\infty,$ for any
$m\in {\Bbb N}$ there exists $i=i(m)$ such that $f(x_i)\in d_m$ for
$i\geqslant i(m).$ Relabeling the sequence $x_i$ if necessary, we
may consider that $f(x_i)\in d_i$ for any $i\in {\Bbb N}.$
Similarly, we may consider that $f(y_i)\in g_i$ for any $i\in {\Bbb
N}.$ Let $\alpha_i:[0, 1]\rightarrow D^{\,\prime}$ be a path joining
$f(x_1)$ and $f(x_i)$ in $d_1,$ and let $\beta_i:[0, 1]\rightarrow
D^{\,\prime}$ be a path joining $f(y_1)$ and $f(y_i)$ in $g_1.$ Let
$\widetilde{\alpha_i}:[0, c_1)\rightarrow D$ and
$\widetilde{\beta_i}:[0, c_2)\rightarrow D^{\,\prime}$ be maximal
$f$-liftings of paths $\alpha_i$ and $\beta_i$ starting at $x_i$ and
$y_i,$ respectively (they exist by Proposition~\ref{pr3}). Due to
the same proposition, only one of the following two situations are
possible:

\medskip
1) $\widetilde{\alpha_i}(t)\rightarrow x_1\in D$ as $t\rightarrow
c_1,$ and $c_1=1$ and $f(\widetilde{\alpha_i}(1))=f(x_1),$ or

\medskip
2) $\widetilde{\alpha_i}(t)\rightarrow \partial D$ as $t\rightarrow
c_1.$

\medskip
Let us to show that the case 2) is impossible. We reason by
contradiction: let $\widetilde{\alpha_i}(t)\rightarrow \partial D$
as $t\rightarrow c_1.$ Let us choose an arbitrary sequence $t_m\in
[0, c_1)$ such that $t_m\rightarrow c_1-0$ as $m\rightarrow\infty.$
Since the space $\overline{{\Bbb R}^n}$ is compact, the boundary
$\partial D$ is also compact as a closed subset of the compact
space. Then there exists $w_m\in \partial D$ such that
\begin{equation}\label{eq7B}
h(\widetilde{\alpha_i}(t_m), \partial
D)=h(\widetilde{\alpha_i}(t_m), w_m) \rightarrow 0\,,\qquad
m\rightarrow \infty\,.
\end{equation}
Due to the compactness of $\partial D$, we may assume that
$w_m\rightarrow w_0\in \partial D$ as $m\rightarrow\infty.$
Therefore, by the relation~(\ref{eq7B}) and by the triangle
inequality
\begin{equation}\label{eq8B}
h(\widetilde{\alpha_i}(t_m), w_0)\leqslant
h(\widetilde{\alpha_i}(t_m), w_m)+h(w_m, w_0)\rightarrow 0\,,\qquad
m\rightarrow \infty\,.
\end{equation}
On the other hand,
\begin{equation}\label{eq9B}
f(\widetilde{\alpha_i}(t_m))=\alpha_i(t_m)\rightarrow \alpha(c_1)
\,,\quad m\rightarrow\infty\,,
\end{equation}
because by the construction the path $\alpha_i(t),$ $t\in [0, 1],$
lies in $D^{\,\prime}\setminus C(f, \partial D)$ together with its
finite ones points. At the same time, by~(\ref{eq8B})
and~(\ref{eq9B}) we have that $\alpha_i(c_1)\in C(f,
\partial D).$ The inclusion $|\alpha_i|\subset D^{\,\prime}\setminus
C(f, \partial D)$ contradicts the relation $\alpha_i(c_1)\in C(f,
\partial D),$ where, as usual, give a path $\gamma:[a,
b]\rightarrow \overline{{\Bbb R}^n}$ we use the notation
$$|\gamma|=\{x\in\overline{{\Bbb R}^n}: \exists\,t\in [a,
b]:\gamma(t)=x\}\,$$ for the {\it locus} of $\gamma.$ The resulting
contradiction indicates the impossibility of the second situation.

\medskip
Therefore, the first situation is fulfilled:
$\widetilde{\alpha_i}(t)\rightarrow x_1\in D$ as
$i\rightarrow\infty,$ and $c_1=1$ and
$f(\widetilde{\alpha_i}(1))=f(x_1).$ In other words, the $f$-lifting
$\widetilde{\alpha_i}$ is complete, i.e., $\widetilde{\alpha_i}:[0,
1]\rightarrow D.$ Similarly, the path $\beta_i$ has a complete
$f$-lifting $\widetilde{\beta_i}:[0, 1]\rightarrow D.$

\medskip
Note that, the points $f(x_1)$ and $f(y_1)$ have only a finite
number of pre-images in $D.$ Let us prove this from the opposite:
suppose that there is a sequence $z_i,$ $i=1,2,\ldots ,$ such that
$f(z_i)=f(x_1).$ Due to the compactness of $\overline{{\Bbb R}^n}$
we may assume that $z_i$ converges to some point $z_0$ as
$i\rightarrow\infty.$ Since $f$ is discrete, $z_0\in
\partial D.$ But then $f(x_1)\in C(f, \partial D),$ which contradicts
the definition of $f(x_1).$ The resulting contradiction indicates a
finite number of pre-images of the point $f(x_1)$ in $D$ under the
mapping $f.$ A similar conclusion may be done with respect to
$f(y_1).$

\medskip
Therefore, there exists $R_0>0$ such that $\widetilde{\alpha_i}(1),
\widetilde{\beta_i}(1)\in D\setminus B(x_0, R_0)$ for all
$i=1,2,\ldots .$ Since the boundary of the domain $D$ is weakly
flat, for every $P>0$ there is $i=i_P\geqslant 1$ such that
\begin{equation}\label{eq7}
M(\Gamma(|\widetilde{\alpha_i}|, |\widetilde{\beta_i}|,
D))>P\qquad\forall\,\,i\geqslant i_P\,.
\end{equation}
Let us to show that the condition~(\ref{eq7}) contradicts the
definition of the mapping $f$ in~(\ref{eq2*A}).

Indeed, let $\gamma\in \Gamma(|\widetilde{\alpha_i}|,
\widetilde{\beta_i}, D),$ then $\gamma:[0, 1]\rightarrow D,$
$\gamma(0)\in |\widetilde{\alpha_i}|$ and $\gamma(1)\in
\widetilde{\beta_i}.$ In particular, $f(\gamma(0))\in |\alpha_i|$
and $f(\gamma(1))\in |\beta_i|.$ In this case, due to~(\ref{eq4})
and~(\ref{eq7}) we obtain that $|f(\gamma)|\cap d_1\ne\varnothing
\ne |f(\gamma)|\cap(D^{\,\prime}\setminus d_1)$ for
$i\geqslant\max\{i_1, i_2\}.$ Due to~\cite[Theorem~1.I.5.46]{Ku$_2$}
$|f(\gamma)|\cap
\partial d_1\ne\varnothing,$ in other words, $|f(\gamma)|\cap S(z_0,
r_1)\ne\varnothing$ because $\partial d_1\cap D^{\,\prime}\subset
\sigma_1\subset S(z_0, r_1)$ by the definition of the cut
$\sigma_1.$ Let $t_1\in (0,1)$ be such that $f(\gamma(t_1))\in
S(z_0, r_1)$ and $f(\gamma)|_1:=f(\gamma)|_{[t_1, 1]}.$ Without loss
of generality, we may consider that $f(\gamma)|_1\subset {\Bbb
R}^n\setminus B(z_0, r_1).$ Arguing similarly for $f(\gamma)|_1,$ we
may find $t_2\in (t_1,1)$ such that $f(\gamma(t_2))\in S(z_0, r_0).$
Set $f(\gamma)|_2:=f(\gamma)|_{[t_1, t_2]}.$ Then the path
$f(\gamma)|_2$ is a subpath of $f(\gamma),$ in addition,
$f(\gamma)|_2\in \Gamma(S(z_0, r_1), S(z_0, r_0), D^{\,\prime}).$
Without loss of generality we may consider that $f(\gamma)|_2\subset
B(z_0, r_0).$ Thus $\Gamma(|\widetilde{\alpha_i}|,
\widetilde{\beta_i}, D)>\Gamma_f(z_0, r_1, r_0).$
By the latter relation, by the minorization of the modulus (see,
e.g., \cite[Theorem~1(c)]{Fu})
\begin{equation}\label{eq10}
\Gamma(|\widetilde{\alpha_i}|, |\widetilde{\beta_i}|, D)
>\Gamma_f(z_0, r_1, r_0)\,.
\end{equation}
In turn, by~(\ref{eq10}) we have that
$$M(\Gamma(|\widetilde{\alpha_i}|, |\widetilde{\beta_i}|, D))\leqslant$$
\begin{equation}\label{eq11}
\leqslant M(\Gamma_f(z_0, r_1, r_0))\leqslant \int\limits_{A}
Q(y)\cdot \eta^n (|y-z_0|)\, dm(y)\,,
\end{equation}
where $A=A(z_0, r_1, r_2)$ and $\eta$ is an arbitrary Lebesgue
measurable function satisfying the relation~(\ref{eqA2}). We use the
following standard conventions: $a/\infty=0$ for $a\ne\infty,$
$a/0=\infty$ for $a>0$ and $0\cdot\infty=0$ (see, e.g.,
\cite[3.I]{Sa}). Set $\widetilde{Q}(y)=\max\{Q(y), 1\},$ and let
\begin{equation}\label{eq13}
I=\int\limits_{r_1}^{r_2}\frac{dt}{t\widetilde{q}_{z_0}^{1/(n-1)}(t)}\,,
\end{equation}
where
\begin{equation}\label{eq12}
\widetilde{q}_{z_0}(r)=\frac{1}{\omega_{n-1}r^{n-1}}\int\limits_{S(z_0,
r)}\widetilde{Q}(y)\,d\mathcal{H}^{n-1}(y)\,, \end{equation}
and $\omega_{n-1}$ is the area of the unit sphere ${\Bbb S}^{n-1}$
in ${\Bbb R}^n.$ By assumption, there is a set $E\subset [r_1 ,
r_2]$ of positive linear Lebesgue measure such that $q_{z_0}(t)$ is
finite at of all $t\in E.$ Consequently, $\widetilde{q}_{z_0}(t)$ is
finite at of all $t\in E,$ as well. Therefore, $I\ne 0$
in~(\ref{eq13}). Besides that, $I\ne\infty,$ because
$$I\leqslant \log\frac{r_2}{r_1}<\infty\,.$$
Put $\eta_0(t)=\frac{1}{Itq_{z_0}^{1/(n-1)}(t)}$ and observe that
$\eta_0$ satisfies the relation~(\ref{eqA2}). Let us substitute this
function in the right part of the inequality~(\ref{eq11}) and apply
Fubini's theorem. We will have that
\begin{equation}\label{eq14}
M(\Gamma(|\widetilde{\alpha_i}|, |\widetilde{\beta_i}|, D))
\leqslant \frac{\omega_{n-1}}{I^{n-1}}<\infty\,.
\end{equation}
The relation~(\ref{eq14}) contradicts~(\ref{eq7}). The contradiction
obtained above shows that the assumption in~(\ref{eq1B}) is wrong.

\medskip
It remains to show the equality
$\overline{f}(\overline{D})=\bigcup\limits_{i=1}^N\overline{D_i}_P.$
Obviously, $\overline{f}(\overline{D})\subset
\bigcup\limits_{i=1}^N\overline{D_i}_P.$ Let us to show the inverse
inclusion. Indeed, let $y_0\in
\bigcup\limits_{i=1}^N\overline{D_i}_P.$ Then either $y_0\in
D^{\,\prime},$ or $y_0\in \bigcup\limits_{i=1}^N E_{D_i}.$ If
$y_0\in D^{\,\prime},$ then $y_0=f(x_0)$ and $y_0\in
\overline{f}(\overline{D}),$ because $f$ maps $D$ onto
$D^{\,\prime}$ by the assumption. Finally, let $y_0\in
\bigcup\limits_{i=1}^NE_{D_i}.$ Now, there exists $1\leqslant
i_0\leqslant N$ such that $y_0\in E_{D_{i_0}}.$ Due to the
regularity of $D_{i_0},$ there exists a sequence $y_k\in D_{i_0}$
such that $\rho_{i_0}(y_k, y_0)\rightarrow 0$ as
$k\rightarrow\infty,$ where $\rho_{i_0}$ is a metric in
$\overline{D_{i_0}}_P$ defined similarly to~(\ref{eq5M}). Obviously,
$D^{\,\prime}$ is a bounded domain. Thus, $y_k$ contains a
convergent subsequence $y_{k_l}\rightarrow y_*\in
\overline{D^{\,\prime}}$ as $l\rightarrow\infty.$ Without loss of
generality we may assume that $y_k\rightarrow y_*$ as
$k\rightarrow\infty.$ Observe that, $y_*\not\in D_{i_0},$ because
$y_0\in E_{D_{i_0}}.$ Thus, $y_*\in C(f, \partial D).$ Consequently,
$y_k=f(x_k),$ $x_k\in D\setminus f^{\,-1}(C(f,
\partial D)\cap D^{\,\prime})$ for $k=1,2,\ldots, $ $x_k\rightarrow
x_0$ as $k\rightarrow \infty.$ Since by the proved above $f$ is
continuous in $D\setminus f^{\,-1}(C(f,
\partial D)\cap D^{\,\prime}),$ we have that
$f(x_0)=y_0\in \overline{f}(\partial D)\subset
\overline{f}(\overline{D}).$ Theorem is proved.~$\Box$

\medskip
{\it Proof of Corollary~\ref{cor1}.} Let $0<r_0:=\sup\limits_{y\in
D^{\,\prime}}|y-y_0|.$ We may assume that $Q$ is extended by zero
outside $D^{\,\prime}.$ By the Fubini theorem (see, e.g.,
\cite[Theorem~8.1.III]{Sa}) we obtain that
$$\int\limits_{r_1<|y-y_0|<r_2}Q(y)\,dm(y)=\int\limits_{r_1}^{r_2}
\int\limits_{S(y_0, r)}Q(y)\,d\mathcal{H}^{n-1}(y)dr<\infty\,.$$
This means the fulfillment of condition~1) in
Theorem~\ref{th3}.~$\Box$

\section{The case of finite connectedness of the domain on the cluster set}

As we see from the formulation of Theorem~\ref{th3}, under the
stated conditions, the continuous extension of mappings
$f:\overline{D}\rightarrow \overline{D^{\,\prime}}_P$ in terms of
the space $\overline{D^{\,\prime}}_P$ is not guaranteed. We are
talking only about extension of mapping $f$ in terms of a special
metric. In order for a more natural continuous extension of the
mapping $f$ to be satisfied, some more special conditions are
needed, which we will now consider.

\medskip
Let $E\subset \overline{D}.$ We say that $D$ is {\it finitely
connected at the point $z_0\in E,$} if for any neighborhood
$\widetilde{U}$ of $z_0$ there is a neighborhood
$\widetilde{V}\subset \widetilde{U}$ of $z_0$ such that $(D\cap
\widetilde{V})\setminus E$ consists of finite number of components.
We say that $D$ is {\it finitely connected on $E,$} if $D$ is
finitely connected at any point $z_0\in E.$ The following theorem is
true.

\medskip
\begin{theorem}\label{th1}
{\it\, Let $D$ and $D^{\,\prime}$ be domains in ${\Bbb R}^n,$
$n\geqslant 2,$ and let $D$ be a domain with a weakly flat boundary.
Suppose that $f$ is open discrete mapping of $D$ onto $D^{\,\prime}$
satisfying the relation~(\ref{eq2*A}) at each point $y_0\in
D^{\,\prime}.$ In addition, assume that the following conditions are
fulfilled:

\medskip
1) for each point $y_0\in \partial D^{\,\prime}$ there is
$0<r_0:=\sup\limits_{y\in D^{\,\prime}}|y-y_0|$ such that for any
$0<r_1<r_2<r_0:=\sup\limits_{y\in D^{\,\prime}}|y-y_0|$ there exists
a set $E\subset[r_1, r_2]$ of positive linear Lebesgue measure such
that $Q$ is integrable on $S(y_0, r)$ for $r\in E;$

\medskip
2) $D^{\,\prime}$ is a regular domain and, in addition,
$D^{\,\prime}\setminus C(f,
\partial D))$ is finitely connected on $C(f,
\partial D)\cap D^{\,\prime},$ i.e.,  for each point $z_0\in C(f, \partial D)\cap
D^{\,\prime}$ and for any neighborhood $U$ of this point there
exists a neighborhood $V\subset U$ of this point such that the set
$V\setminus C(f,
\partial D)$ consists of a finite number of
components;

\medskip
3) the set $f^{\,-1}(C(f, \partial D)\cap D^{\,\prime})$ is nowhere
dense in $D.$

\medskip
Then the mapping $f$ has a continuous extension
$\overline{f}:\overline{D}\rightarrow
D^{\,\prime}\cup\bigcup\limits_{i=1}^N\overline{D_i}_P$ by the
metric $\rho$ defined in~(\ref{eq1}). Moreover,
$\overline{f}(\overline{D})=D^{\,\prime}\cup\bigcup\limits_{i=1}^N\overline{D_i}_P.$
}
\end{theorem}

\medskip
\begin{proof}
We fix $x_0\in\partial D.$ Let us show that the mapping $f$ has a
continuous extension to a point $x_0.$ Using the M\"{o}bius
transformation $\varphi:\infty\mapsto 0$ if necessary, due to the
invariance of the modulus $M$ in~(\ref{eq2*A})
(see~\cite[Theorem~8.1]{Va}), we may assume that $x_0\ne\infty.$

\medskip
Let us carry out the proof from the opposite, namely, suppose that
$f$ has no a continuous extension to the point $x_0.$ Due to the
compactness of the space $\overline{D^{\,\prime}}_P,$ there exist
sequences $x_i, y_i\in D,$ $i=1,2,\ldots ,$ such that $x_i,
y_i\rightarrow x_0$ as $i\rightarrow\infty,$ and
\begin{equation}\label{eq1C}
\rho(f(x_i), f(y_i))\geqslant a>0
\end{equation}
for some $a>0$ and all $i\in {\Bbb N},$ where $\rho$ is a metric
in~(\ref{eq5M}). Again, by the compactness of
$\overline{D^{\,\prime}}_P,$ we may assume that the sequences
$f(x_i)$ and $f(y_i)$ converge as $i\rightarrow\infty$ to $P_1$ and
$P_2\in \overline{D^{\,\prime}}_P,$ respectively.

\medskip
First of all, we note that the points $x_i$ and $y_i,$
$i=1,2,\ldots, $ may be chosen such that $x_i, y_i\not\in
f^{\,-1}(C(f,
\partial D)\cap D^{\,\prime}).$ Indeed, since under condition~3) the set
$f^{\,-1}(C(f, \partial D)\cap D^{\,\prime})$ is nowhere dense in
$D,$ there exists a sequence $x_{ki}\in D\setminus(f^{\,-1}(C(f,
\partial D)\cap D^{\,\prime})),$
$k=1,2,\ldots ,$ such that $x_{ki}\rightarrow x_i$ as
$k\rightarrow\infty.$ Put $\varepsilon>0.$ Due to the continuity of
the mapping $f$ at the point $x_i,$ for the number $i\in {\Bbb N}$
there is a number $k_i\in {\Bbb N}$ such that $h(f(x_{k_ii}),
f(x_i))<\frac{1}{2^i}.$ So, by the triangle inequality
$$h(f(x_{k_ii}), z_1)\leqslant h(f(x_{k_ii}), f(x_i))+
h(f(x_i), z_1)\leqslant \frac{1}{2^i}+\varepsilon\,,$$
$i\geqslant i_0=i_0(\varepsilon),$ since $f(x_i)\rightarrow z_1$ as
$i\rightarrow\infty$ and by the choice of $x_i$ and $z_1.$
Therefore, $x_i\in D$ may be replaced by $x_{k_ii}\in
D\setminus(f^{\,-1}(C(f,
\partial D)\cap D^{\,\prime})),$ as required. We may reason similarly for $y_i.$

\medskip
Let us consider two possible cases:

\medskip
{\bf 1). At least one of the sequences $f(x_i)$ of $f(y_i)$ converge
to some boundary element $P_1\in E_{D^{\,\prime}}$ as
$i\rightarrow\infty.$} Let, for example, $f(x_i)\rightarrow P_1\in
E_{D^{\,\prime}}$ as $i\rightarrow\infty.$

\medskip
Let $f(y_i)\rightarrow P_1\in {D^{\,\prime}}_P$ as
$i\rightarrow\infty.$ Let $\sigma_i$ and $\sigma^{\,\prime}_i,$
$i=0,1,2,\ldots, $ be sequences of cuts corresponding $P_1$ and
$P_2.$ In the case where $P_2$ is an inner point of $D^{\,\prime}$
we will consider the sequence of spheres $\sigma^{\,\prime}_i=S(P_2,
R_i)$ as the corresponding system of cuts, where $R_i>0,$
$R_i\rightarrow 0$ as $i\rightarrow\infty.$

\medskip
We may consider that the cuts $\sigma_i,$ $i=0,1,2,\ldots, $ lie on
the spheres $S(z_0, r_i)$ centered at some point $z_0\in
\partial D^{\,\prime},$ where $r_i\rightarrow 0$ as $i\rightarrow\infty$
(such a sequence $\sigma_i$ exists due to~\cite[Lemma~3.1]{IS},
cf.~\cite[Lemma~1]{KR}). Let $d_i$ be a sequence of nested domains
in $D^{\,\prime},$ corresponding to the sequence of cuts $\sigma_i.$
If $P_2\in E_{D^{\,\prime}},$ we consider the sequence of nested
domains $g_i,$ $i=0,1,2,\ldots, $ in $D^{\,\prime},$ corresponding
to $\sigma^{\,\prime}_i.$ Obviously, we may consider that $d_0\cap
g_0=\varnothing\,.$ Otherwise, if $P_2\in D,$ we may choose $R_*>0$
such that $B(P_2, R_*)\cap d_0=\varnothing.$ By the conditions of
Theorem~\ref{th1}, there exists a neighborhood $V\subset B(P_2,
R_*)$ such that $V\setminus C(f, D)$ consists a finite number of
components. Thus, at least one a component $K_2$ contains infinitely
many elements $f(y_i),$ $i=1,2,\ldots. $ Without loss of generality,
relabeling $f(y_i)$ if necessary, we may assume that $f(y_i)\in K_2$
for any $i\in {\Bbb N}.$ Set $g_0:=K_2.$ Thus, in any case $P_2\in
E_{D^{\,\prime}}$ or $P_2\in D^{\,\prime},$
\begin{equation}\label{eq4A}
d_0\cap g_0=\varnothing\,.
\end{equation}
Since $f(x_i)$ converges to $P_1$ as $i\rightarrow\infty,$ for any
$m\in {\Bbb N}$ there exists $i=i(m)$ such that $f(x_i)\in d_m$ for
$i\geqslant i(m).$ Relabeling the sequence $x_i$ if necessary, we
may consider that $f(x_i)\in d_i$ for any $i\in {\Bbb N}.$ We may
consider that $f(y_i)\in g_0$ for any $i\in {\Bbb N}.$

Since by the definition of the prime end,
$\bigcap\limits_{k=1}^{\infty}d_k=\varnothing,$ way consider that
$f(x_i)\not \in d_0$ for $i\in {\Bbb N}.$ Let $\alpha_i:[0,
1]\rightarrow D^{\,\prime}$ be a path joining $f(x_1)$ and $f(x_i)$
in $d_1,$ and let $\beta_i:[0, 1]\rightarrow D^{\,\prime}$ be a path
joining $f(y_1)$ and $f(y_i)$ in $g_0.$ Let
$\widetilde{\alpha_i}:[0, c_1)\rightarrow D$ and
$\widetilde{\beta_i}:[0, c_2)\rightarrow D^{\,\prime}$ be maximal
$f$-liftings of paths $\alpha_i$ and $\beta_i$ starting at $x_i$ and
$y_i,$ respectively (they exist by Proposition~\ref{pr3}). Due to
the same proposition, only one of the following two situations are
possible:

\medskip
1) $\widetilde{\alpha_i}(t)\rightarrow x_1\in D$ as $t\rightarrow
c_1,$ and $c_1=1$ and $f(\widetilde{\alpha_i}(1))=f(x_1),$ or

\medskip
2) $\widetilde{\alpha_i}(t)\rightarrow \partial D$ as $t\rightarrow
c_1.$

\medskip
Arguing similarly to the proof of Theorem~\ref{th3}, we may prove
that the case 2) is impossible.

\medskip
Therefore, the first situation is fulfilled:
$\widetilde{\alpha_i}(t)\rightarrow x_1\in D$ as
$i\rightarrow\infty,$ and $c_1=1$ and
$f(\widetilde{\alpha_i}(1))=f(x_1).$ In other words, the $f$-lifting
$\widetilde{\alpha_i}$ is complete, i.e., $\widetilde{\alpha_i}:[0,
1]\rightarrow D.$ Similarly, the path $\beta_i$ has a complete
$f$-lifting $\widetilde{\beta_i}:[0, 1]\rightarrow D.$

\medskip
We may show that, the points $f(x_1)$ and $f(y_1)$ have only a
finite number of pre-images in $D$ (see the proof of
Theorem~\ref{th3}). Therefore, there exists $R_0>0$ such that
$\widetilde{\alpha_i}(1), \widetilde{\beta_i}(1)\in D\setminus
B(x_0, R_0)$ for all $i=1,2,\ldots .$ Since the boundary of the
domain $D$ is weakly flat, for every $P>0$ there is $i=i_P\geqslant
1$ such that
\begin{equation}\label{eq7C}
M(\Gamma(|\widetilde{\alpha_i}|, |\widetilde{\beta_i}|,
D))>P\qquad\forall\,\,i\geqslant i_P\,.
\end{equation}
Let us to show that the condition~(\ref{eq7C}) contradicts the
definition of the mapping $f$ in~(\ref{eq2*A}). Indeed, arguing
similarly to the proof of Theorem~\ref{th3}, we obtain that
\begin{equation}\label{eq10A}
\Gamma(|\widetilde{\alpha_i}|, |\widetilde{\beta_i}|, D)
>\Gamma_f(z_0, r_1, r_0)\,.
\end{equation}
Now, by~(\ref{eq10A}) we have that
$$M(\Gamma(|\widetilde{\alpha_i}|, |\widetilde{\beta_i}|, D))\leqslant$$
\begin{equation}\label{eq11A}
\leqslant M(\Gamma_f(z_0, r_1, r_0))\leqslant \int\limits_{A}
Q(y)\cdot \eta^n (|y-z_0|)\, dm(y)\,,
\end{equation}
where $A=A(z_0, r_1, r_2)$ and $\eta$ is an arbitrary Lebesgue
measurable function satisfying the relation~(\ref{eqA2}). Set
$\widetilde{Q}(y)=\max\{Q(y), 1\},$ and let
\begin{equation}\label{eq13A}
I=\int\limits_{r_1}^{r_2}\frac{dt}{t\widetilde{q}_{z_0}^{1/(n-1)}(t)}\,,
\end{equation}
where
\begin{equation}\label{eq12A}
\widetilde{q}_{z_0}(r)=\frac{1}{\omega_{n-1}r^{n-1}}\int\limits_{S(z_0,
r)}\widetilde{Q}(y)\,d\mathcal{H}^{n-1}(y)\,, \end{equation}
and $\omega_{n-1}$ is the area of the unit sphere ${\Bbb S}^{n-1}$
in ${\Bbb R}^n.$ By assumption, there is a set $E\subset [r_1 ,
r_2]$ of positive linear Lebesgue measure such that $q_{z_0}(t)$ is
finite at of all $t\in E.$ Consequently, $\widetilde{q}_{z_0}(t)$ is
finite at of all $t\in E,$ as well. Therefore, $I\ne 0$
in~(\ref{eq13}). Besides that, $I\ne\infty,$ because
$$I\leqslant \log\frac{r_2}{r_1}<\infty\,.$$
Put $\eta_0(t)=\frac{1}{Itq_{z_0}^{1/(n-1)}(t)}$ and observe that
$\eta_0$ satisfies the relation~(\ref{eqA2}). Let us substitute this
function in the right part of the inequality~(\ref{eq11}) and apply
Fubini's theorem. We will have that
\begin{equation}\label{eq14A}
M(\Gamma(|\widetilde{\alpha_i}|, |\widetilde{\beta_i}|, D))
\leqslant \frac{\omega_{n-1}}{I^{n-1}}<\infty\,.
\end{equation}
The relation~(\ref{eq14A}) contradicts~(\ref{eq7C}). The
contradiction obtained above shows that the assumption
in~(\ref{eq1C}) is wrong.

\medskip
{\bf 2). All of the sequences $f(x_i)$ of $f(y_i)$ converge to some
inner elements $z_1$ and $z_2\in D^{\,\prime}$ as
$i\rightarrow\infty.$} Here we may assume that $z_1\ne\infty.$

\medskip
By the condition~2), given neighborhoods $U_1$ and $U_2$ of the
points $z_1$ and $z_2\in C(f,
\partial D),$ there are neighborhoods $V_1\subset U_1$
and $V_2\subset U_2$ of these same points such that each of the sets
$W_1:=V_1\setminus C(f,
\partial D)$ and $W_2:=V_2\setminus C(f,
\partial D)$ consist of a finite number of components. Observe that, we may choose
open $V_1$ and $V_2.$ Then $W_1$ and $W_2$ are also open as the
intersection of two open sets. Therefore, each component of $W_1$
and $W_2$ is open (see~\cite[Theorem~4.II.49.6]{Ku$_2$}) and path
connected (see, e.g.,~\cite[Proposition~13.2]{MRSY}).

\medskip
Without loss of generalization, we may assume that the number
$r_0=r_0(z_1)>0$ numbers $0<r_1<r_2<r_0$ from the condition of the
theorem are such that
\begin{equation}\label{eq8}
U_1\subset B(z_1, r_1),\qquad \overline{B(z_1, r_2)}\cap
\overline{U_2}=\varnothing\,.
\end{equation}
Passing to subsequences if necessary and taking into account that
$W_1$ and $W_2$ have a finite number of components, we may assume
that $f(x_i),$ $i=1,2,\ldots ,$ belong to some single component
$K_1$ of the set $W_1.$ Similarly, we may assume that $f(y_i),$
$i=1,2,\ldots ,$ belong to single component $K_2$ of the set $W_2.$
Let us to join the points $f(x_i)$ and $f(x_1)$ by a path
$\alpha_i:[0, 1]\rightarrow D^{\,\prime},$ and the points $f(y_i)$
and $f(y_1)$ by a path $\beta_i:[0, 1]\rightarrow D^{\,\prime}$ so
that $|\alpha_i|\subset K_1$ and $|\beta_i|\subset K_2$ for
$i=1,2,\ldots .$ Let $\widetilde{\alpha_i}:[0, c_1)\rightarrow D$
and $\widetilde{\beta_i}:[0, c_2)\rightarrow D^{\,\prime}$ be
maximal $f$-liftings of paths $\alpha_i$ and $\beta_i$ starting at
$x_i$ and $y_i,$, respectively (they exist by
Proposition~\ref{pr3}). Due to the same proposition, only one of the
following two situations are possible:

\medskip
1) $\widetilde{\alpha_i}(t)\rightarrow x_1\in D$ as $t\rightarrow
c_1,$ and $c_1=1$ and $f(\widetilde{\alpha_i}(1))=f(x_1),$ or

\medskip
2) $\widetilde{\alpha_i}(t)\rightarrow \partial D$ as $t\rightarrow
c_1.$

\medskip
Arguing similarly to the proof of Theorem~\ref{th3}, we may prove
that the case 2) is impossible.
Thus,$\widetilde{\alpha_i}(t)\rightarrow x_1\in D$ as
$i\rightarrow\infty,$ and $c_1=1$ and
$f(\widetilde{\alpha_i}(1))=f(x_1).$ In other words, the $f$-lifting
$\widetilde{\alpha_i}$ is complete, i.e., $\widetilde{\alpha_i}:[0,
1]\rightarrow D.$ Similarly, the path $\beta_i$ has a complete
$f$-lifting $\widetilde{\beta_i}:[0, 1]\rightarrow D.$

\medskip
We may show that, the points $f(x_1)$ and $f(y_1)$ have only a
finite number of pre-images in $D$ (see the proof of
Theorem~\ref{th3}). Therefore, there exists $R_0>0$ such that
$\widetilde{\alpha_i}(1), \widetilde{\beta_i}(1)\in D\setminus
B(x_0, R_0)$ for all $i=1,2,\ldots .$ Since the boundary of the
domain $D$ is weakly flat, for every $P>0$ there is $i=i_P\geqslant
1$ such that
\begin{equation}\label{eq7D}
M(\Gamma(|\widetilde{\alpha_i}|, |\widetilde{\beta_i}|,
D))>P\qquad\forall\,\,i\geqslant i_P\,.
\end{equation}
Let us to show that the condition~(\ref{eq7}) contradicts the
definition of the mapping $f$ in~(\ref{eq2*A}). Indeed,
by~(\ref{eq8}) and by~\cite[Theorem~1.I.5.46]{Ku$_2$}
\begin{equation}\label{eq9C}
f(\Gamma(|\widetilde{\alpha_i}|, |\widetilde{\beta_i}|,
D))>\Gamma(S(z_1, r_1), S(z_1, r_2), A(z_1, r_1, r_2))\,.
\end{equation}
Therefore, it follows from~(\ref{eq9C}) that
\begin{equation}\label{eq10B}
\Gamma(|\widetilde{\alpha_i}|, |\widetilde{\beta_i}|, D)
>\Gamma_f(z_1, r_1, r_2)\,.
\end{equation}
In turn, by~(\ref{eq10B}) we have that
$$M(\Gamma(|\widetilde{\alpha_i}|, |\widetilde{\beta_i}|, D))\leqslant$$
\begin{equation}\label{eq11B}
\leqslant M(\Gamma_f(z_1, r_1, r_2))\leqslant \int\limits_{A}
Q(y)\cdot \eta^n (|y-z_1|)\, dm(y)\,,
\end{equation}
where $A=A(z_1, r_1, r_2)$ and $\eta$ is an arbitrary Lebesgue
measurable function satisfying the relation~(\ref{eqA2}). We use the
following standard conventions: $a/\infty=0$ for $a\ne\infty,$
$a/0=\infty$ for $a>0$ and $0\cdot\infty=0$ (see, e.g.,
\cite[3.I]{Sa}). Set $\widetilde{Q}(y)=\max\{Q(y), 1\},$ and let
\begin{equation}\label{eq13B}
I=\int\limits_{r_1}^{r_2}\frac{dt}{t\widetilde{q}_{z_1}^{1/(n-1)}(t)}\,,
\end{equation}
where
\begin{equation}\label{eq12B}
\widetilde{q}_{z_1}(r)=\frac{1}{\omega_{n-1}r^{n-1}}\int\limits_{S(z_1,
r)}\widetilde{Q}(y)\,d\mathcal{H}^{n-1}(y)\,, \end{equation}
and $\omega_{n-1}$ is the area of the unit sphere ${\Bbb S}^{n-1}$
in ${\Bbb R}^n.$ By assumption, there is a set $E\subset [r_1 ,
r_2]$ of positive linear Lebesgue measure such that $q_{z_1}(t)$ is
finite at of all $t\in E.$ Consequently, $\widetilde{q}_{z_1}(t)$ is
finite at of all $t\in E,$ as well. Therefore, $I\ne 0$
in~(\ref{eq13}). Besides that, $I\ne\infty,$ because
$$I\leqslant \log\frac{r_2}{r_1}<\infty\,.$$
Put $\eta_0(t)=\frac{1}{Itq_{z_1}^{1/(n-1)}(t)}$ and observe that
$\eta_0$ satisfies the relation~(\ref{eqA2}). Let us substitute this
function in the right part of the inequality~(\ref{eq11}) and apply
Fubini's theorem. We have that
\begin{equation}\label{eq14B}
M(\Gamma(|\widetilde{\alpha_i}|, |\widetilde{\beta_i}|, D))
\leqslant \frac{\omega_{n-1}}{I^{n-1}}<\infty\,.
\end{equation}
The relation~(\ref{eq14B}) contradicts~(\ref{eq7D}). The
contradiction obtained above shows that the assumption
in~(\ref{eq1C}) is wrong.

\medskip
It remains to show the equality
$\overline{f}(\overline{D})=D^{\,\prime}_P.$ Obviously,
$\overline{f}(\overline{D})\subset D^{\,\prime}_P.$ Let us to show
the inverse inclusion. Indeed, let $y_0\in D^{\,\prime}_P.$ Then
either $y_0\in D^{\,\prime},$ or $y_0\in E_{D^{\,\prime}}.$ If
$y_0\in D^{\,\prime},$ then $y_0=f(x_0)$ and $y_0\in
\overline{f}(\overline{D}),$ because $f$ maps $D$ onto
$D^{\,\prime}$ by the assumption. Finally, let $y_0\in
E_{D^{\,\prime}}.$ Due to the regularity of $D^{\,\prime}_{P},$
there exists a sequence $y_k\in D^{\,\prime}$ such that $\rho(y_k,
y_0)\rightarrow 0$ as $k\rightarrow\infty,$ $y_k=f(x_k)$ and $x_k\in
D,$ where $\rho$ is a metric in $\overline{D^{\,\prime}}_P$ defined
by~(\ref{eq5M}). Due to the compactness of the space
$\overline{{\Bbb R}^n}$ we may consider that $x_k\rightarrow x_0,$
as $x_0\in\overline{D}.$ Observe that, $x_0\in
\partial D,$ because $f$ is open. Now,
$f(x_0)=y_0\in \overline{f}(\partial D)\subset
\overline{f}(\overline{D}).$ Theorem is proved.~$\Box$
\end{proof}

\medskip
{\bf \noindent Evgeny Sevost'yanov} \\
{\bf 1.} Zhytomyr Ivan Franko State University,  \\
40 Velyka Berdychivs'ka Str., 10 008  Zhytomyr, UKRAINE \\
{\bf 2.} Institute of Applied Mathematics and Mechanics\\
of NAS of Ukraine, \\
19 Henerala Batyuka Str., 84 116 Slov'yansk,  UKRAINE\\
esevostyanov2009@gmail.com

\medskip
{\bf \noindent Victoria Desyatka} \\
Zhytomyr Ivan Franko State University,  \\
40 Velyka Berdychivs'ka Str., 10 008  Zhytomyr, UKRAINE \\
victoriazehrer@gmail.com

\medskip
\noindent{{\bf Zarina Kovba} \\
Zhytomyr Ivan Franko State University,  \\
40 Velyka Berdychivs'ka Str., 10 008  Zhytomyr, UKRAINE \\
e-mail: mazhydova@gmail.com  }

\end{document}